\documentclass{article}
\usepackage[english,russian,ukrainian]{babel} 
\usepackage{latexsym,amssymb,amsmath}
\usepackage[cp1251]{inputenc}
\usepackage{amsfonts,amssymb}
\usepackage{graphicx,graphics,hhline}
\usepackage{euscript}
\begin{document}
\Large
\hspace*{11cm} {\LARGE \it Математика}\\
\noindent УДК 514.517+517.55 \vspace{3mm}\\
\noindent{\textbf{М.~В.~Стефанчук}} \vspace{3mm} \\
\noindent{\bf Екстремальні елементи в
гіперкомплексному просторі}\\

\noindent {\it (Представлено членом-кореспондентом НАН України
Ю.~Ю.~Трохимчуком)}

\vspace{5mm} У роботі досліджуються екстремальні елементи та
\emph{h}-оболонка множин в \emph{n}-вимірному гіперкомплексному
просторі $\mathbb{H^{\emph{n}}}$. Вводиться клас
$\mathbb{H}$-квазіопуклих множин, які включають в себе сильно
гіперкомплексно опуклі множини та є замкненими відносно перетинів.

\textbf{Ключові слова:} гіперкомплексно опукла множина, сильно
гіперкомплексно опукла множина, \emph{h}-оболонка множини,
\emph{h}-екстремальна точка, \emph{h}-екстремальний промінь,
$\mathbb{H}$-квазіопукла множина.

\vspace{5mm}

Дана робота присвячена дослідженню екстремальних елементів та
\emph{h}-оболонки множини в гіперкомплексному просторі, а також
побудові класу гіперкомплексно опуклих множин, які включають в себе
сильно гіперкомплексно опуклі множини та є замкненими відносно
перетинів. Такі множини будемо називати $\mathbb{H}$-квазіопуклими.

Будемо розглядати \emph{n}-вимірний гіперкомплексний простір
$\mathbb{H}^{n}$, $n\in\mathbb{N}$, який є прямим добутком
\emph{n}-копій тіла кватерніонів $\mathbb{H}$. Нехай
$E\subset\mathbb{H}$~-- довільна множина, яка містить початок
координат $O=\{0,0,...,0\}$. Покладемо $x=(x_{1},x_{2},...,x_{n})$,
$h=(h_{1},h_{2},...,h_{n})$, $\langle
x,h\rangle=x_{1}h_{1}+x_{2}h_{2}+...+x_{n}h_{n}$. Множина
$E^{*}=\{h|\langle h,x\rangle\neq1 \forall x\in E\}$ називається
\emph{спряженою} множиною до множини \emph{Е}. [7]

\textbf{Означення 1.} Множина $E\subset\mathbb{H}^{n}$ називається
\emph{гіперкомплексно опуклою}, якщо для довільної точки
$x_{0}\in\mathbb{H}^{n}\setminus E$ існує гіперплощина, яка
проходить через точку $x_{0}$ і не перетинає \emph{Е}. [7]

Нагадаємо, що оскільки в алгебрі кватерніонів множення
некомутативне, то надалі розглядатимемо праві гіперплощини, тобто
точку зі змінною координатою множимо на фіксовану точку справа.

\textbf{Означення 2.} Множина $E\subset\mathbb{H}^{n}$ називається
\emph{сильно гіперкомплексно опуклою}, якщо довільний її перетин
гіперкомплексною прямою $\gamma$ ациклічний, тобто
$\widetilde{H}^{i}(\gamma\bigcap E)=0, \forall i$, де
$\widetilde{H}^{i}(\gamma\bigcap E)$~-- зведена група когомологій
Александрова-Чеха множини $\gamma\bigcap E$ з коефіцієнтами в групі
цілих чисел. [7]

У [8] доведено, що сильно гіперкомплексно опуклі компакти будуть
гіперкомплексно опуклими.

Нехай $E\subset\mathbb{H}$~-- довільна множина. Доповнення до
об'єднання необмежених компонент множини $\mathbb{H}\setminus E$
називається \emph{h-комбінацією точок множини} \emph{Е} та
позначається $[E]$. Якщо \emph{Е}~-- довільна множина в просторі
$\mathbb{H}^{n},n>1$, то скажемо, що точка \emph{х} належить
\emph{h}-комбінації точок з \emph{Е}, якщо існує перетин множини
\emph{Е} гіперкомплексною прямою $\gamma$ такий, що
$x\in[E\bigcap\gamma]$. Множину таких точок з $\mathbb{H}^{n}$
називають \emph{h}\emph{-комбінацією точок} \emph{Е} і позначають
$[E]$; \emph{m}-кратну \emph{h}-комбінацію визначають за індукцією
$[E]^{m}=[[E]^{m-1}]$. [3]

\textbf{Означення 3.} \emph{h}\emph{-оболонкою множини}
$E\subset\mathbb{H}^{n}$ називається множина
$\widehat{E}=\bigcap\limits_{\pi}\pi^{-1}[\pi(E)]$, де
$\pi:\mathbb{H}^{n}\rightarrow\lambda$~-- всеможливі лінійні
проекції множини на гіперкомплексні прямі, $[\pi(E)]$~--
\emph{h}-комбінація точок множини $\pi(E)$, а
$\pi^{-1}[\pi(E)]=\{x\in\mathbb{H}^{n}|\pi(x)\in\pi(E)\}$~-- її
повний прообраз. [3], [7]

Наступна теорема стверджує, що для довільної множини множина точок
її \emph{h}-оболонки співпадає з \emph{h}-комбінацією точок цієї
множини.

\textbf{Теорема 1.} Якщо множина \emph{Е} є \emph{h}-оболонкою, то
$E=[E]$.

\emph{Доведення.} Нехай $x\in[\lambda\bigcap E]$ для деякої
гіперкомплексної прямої $\lambda$. Тоді очевидно буде виконуватися
включення $\pi(x)\in[\pi(\lambda\bigcap E)]$ для всіх проекцій
$\pi$, тому що обмеження будь-якої проекції $\pi$ на кожну пряму є
або гомеоморфізмом, або проекцією в точку.

Дослідимо більш детально процедуру утворення \emph{h}-оболонки
множини \emph{Е}. За властивістю 19 гіперкомплексно опуклих множин
[7] множина $[\pi(E)]$ гомеоморфна доповненню в
$\mathbb{H}^{o}=\lambda^{o}$ ($\mathbb{H}^{o}$~-- одноточкова
компактифікація прямої) до зв'язної компоненти, що містить початок
координат, деякого перерізу $\lambda\bigcap E^{*}$, де $\lambda$~--
гіперкомплексна пряма, а $E^{*}$~-- множина, спряжена до множини
\emph{Е}. Позначимо цю компоненту $]\lambda\bigcap E^{*}[$. Крім
цього дану компоненту можна зобразити ще й так
$(\pi^{-1}[\pi(E)])^{*}$. Взагалі кажучи, об'єднання
$\bigcup\limits_{\lambda}]\lambda\bigcap E^{*}[$ може не бути
гіперкомплексно опуклою множиною. Однак за властивостями 9 та 15
гіперкомплексно опуклих множин [7] справедливі наступні рівності
$$(\bigcup\limits_{\pi}(\pi^{-1}[\pi(E)])^{*})^{*}=\bigcap\limits_{\pi}(\pi^{-1}[\pi(E)])^{**}=\bigcap\limits_{\pi}\pi^{-1}[\pi(E)]=\widehat{E}.$$
Звідси випливає наступна теорема, яка дає ще один спосіб побудови
\emph{h}-оболонки.

\textbf{Теорема 2.} Для довільної множини $E\subset\mathbb{H}^{n}$
її \emph{h}-оболонку можна зобразити у вигляді
$\widehat{E}=(\bigcup\limits_{\lambda}]\lambda\bigcap E^{*}[)^{*}$.

Доведення випливає з виконання рівності
$$\widehat{E}=(\bigcup\limits_{\pi}(\pi^{-1}[\pi(E)])^{*})^{*}=(\bigcup\limits_{\lambda}]\lambda\bigcap E^{*}[)^{*}.$$

Варто відмітити, що спряжена множина до \emph{h}-оболонки не завжди
буде \emph{h}-оболонкою. Розглянемо наступний приклад.

\textbf{Приклад 1.} Нехай \emph{К}~-- компакт, який лежить у
тривимірному евклідовому просторі
$\mathbb{R}^{3}=\{(y_{1},z_{1},y_{2})|x_{1}=y_{1}+iz_{1}+ju_{1}+kt_{1},x_{2}=y_{2}+iz_{2}+ju_{2}+kt_{2},(x_{1},x_{2})\in\mathbb{H}^{2}\}$,
$K=\{(y_{1},z_{1},y_{2})|(y_{1}-|z_{1}|=0,0\le
y_{1}\le1)\wedge[[(y_{2}=y_{1})\vee(y_{2}=-2y_{1}+3),z_{1}\le0]\vee[(y_{2}=2y_{1})\vee(y_{2}=-y_{1}+3),z_{1}>0]]\}$.

Тоді при канонічній проекції $\pi_{1}(K)$ на пряму $x_{2}=0$
отримаємо ациклічний компакт
$K_{1}=\{(y_{1},z_{1},0,0)|y_{1}-|z_{1}|=0,0\le y_{1}\le1\}$, а при
проекції $\pi_{2}(K)$ на пряму $x_{2}$, паралельно прямій, яка
проходить через точки $(1-i,1)$, $(1+i,2)$, отримаємо компакт
$\pi_{1}(K)=K_{2}$, який складається з двох відрізків, оскільки
відрізки $[y_{1}=-z_{1},0\le
y_{1}\le1,y_{2}=y_{1}]\cap\mathbb{R}^{3}$ і $[y_{1}=z_{1},0\le
y_{1}\le1,y_{2}=2y_{1}]\cap\mathbb{R}^{3}$, $[y_{1}=-z_{1},0\le
y_{1}\le1,y_{2}=-2y_{1}+3]\cap\mathbb{R}^{3}$ і $[y_{1}=z_{1},0\le
y_{1}\le1,y_{2}=-y_{1}+3]\cap\mathbb{R}^{3}$ попарно ототожнюються
при цій проекції. Тому $[K_{1}]=K_{1}$ і $[K_{2}]=K_{2}$. Легко
бачити, що $K=\pi_{1}^{-1}(K_{1})\cap\pi_{2}^{-1}(K_{2})$. При всіх
інших проекціях образ $\pi(K)$ буде ненульовим одновимірним циклом.
Тому $\pi(K)\neq[\pi(K)]$ при $\pi\neq\pi_{1}$ або $\pi\neq\pi_{2}$.
Але
$\widehat{K}=\bigcap\limits_{\pi}\pi^{-1}[\pi(E)]\subset\pi_{1}^{-1}(K_{1})\cap\pi_{2}^{-1}(K_{2})=K$.
Отже, $\widehat{K}=K$. Однак, взагалі кажучи,
$\pi(\widehat{K})\neq[\pi(K)]$.

За теоремою 2 $\widehat{K}=(\bigcup\limits_{\gamma}]\gamma\bigcap
K^{*}[)^{*}=(K_{3})^{*}$, де
$K_{3}=\{(y_{1},z_{1},y_{2})|y_{1}-|z_{1}|=0,0\le y_{1}\le
1,(y_{2}=y_{1})\vee(y_{2}=-2y_{1}+3)\}$~-- частина компакта
\emph{К}, що складається з двох відрізків.

З леми 1, яка є аналогом теореми Каратеодорі, враховуючи, що жодна
\emph{h}-екстремальна точка не може бути отримана
\emph{h}-комбінацією інших точок компакта \emph{К}, випливають
наступні наслідки.

\textbf{Лема 1.} \emph{h}-оболонка гіперкомплексно опуклого компакта
в $\mathbb{H}^{n}$ співпадає з сукупністю не більше, ніж
\emph{n}-кратних комбінацій своїх \emph{h}-екстремальних точок. [3]

\textbf{Наслідок 1.} Якщо $K\subset\mathbb{H}^{n}$~-- компакт і
$\widehat{K}$~-- його \emph{h}-оболонка, яка співпадає з
\emph{h}-комбінацією $[K]^{m}$, то кожна \emph{h}-екстремальна точка
множини $\widehat{K}$ належить \emph{К}.

\textbf{Наслідок 2.} Довільну точку простору $\mathbb{H}^{n}$, яка
належить \emph{h}-оболонці $\widehat{K}=[K]^{m}$, можна зобразити у
вигляді не більше, ніж \emph{n}-кратної комбінації точок компакта
\emph{К}.

\textbf{Наслідок 3.} Для того, щоб \emph{h}-оболонка гіперкомплексно
опуклого компакта \emph{К} співпадала з ним, необхідно, щоб всі
перерізи його гіперкомплексними прямими $\gamma$ не містили
тривимірних коциклів, тобто $H^{3}(K\bigcap\gamma)=0$.

\textbf{Наслідок 4.} Для того, щоб \emph{h}-оболонка гіперкомплексно
опуклого компакта \emph{К} співпадала з ним, необхідно, щоб усі
проекції його спряженої множини $K^{*}$ на гіперкомплексні прямі
були зв'язними.

Останній наслідок отримується з попереднього, використовуючи
властивості 19 та 20 гіперкомплексно опуклих множин [7].

\textbf{Приклад 2.} Нехай
$S^{4}_{+}\subset\mathbb{R}^{5}\subset\mathbb{H}^{2}$~-- півсфера у
п'ятивимірному дійсному евклідовому просторі
$\mathbb{R}^{5}=\{(y_{1},z_{1},u_{1},t_{1},y_{2})|x_{1}=y_{1}+iz_{1}+ju_{1}+kt_{1},x_{2}=y_{2}+iz_{2}+ju_{2}+kt_{2},(x_{1},x_{2})\in\mathbb{H}^{2}\}$,
$S^{4}_{+}=\{(y_{1},z_{1},u_{1},t_{1},y_{2})|y_{1}^{2}+z_{1}^{2}+u_{1}^{2}+t_{1}^{2}+y_{2}^{2}=1,z_{1}\ge0\}$.
Легко бачити, що
$(S^{4}_{+})^{**}=K_{+}=\{(y_{1},z_{1},u_{1},t_{1},y_{2})|y_{1}^{2}+z_{1}^{2}+u_{1}^{2}+t_{1}^{2}+y_{2}^{2}\le1,z_{1}\ge0\}$
і $K_{+}$~-- опуклий компакт, який співпадає зі своєю
\emph{h}-оболонкою. Довільний перетин $S^{4}_{+}$ гіперкомплексною
прямою співпадає з перетином $S^{4}_{+}$ дійсною прямою або дійсною
чотирьохвимірною площиною виду $y_{2}=\texttt{const}$ в
$\mathbb{R}^{5}$. Всі такі перетини не містять тривимірних циклів.
Тому $[S^{4}_{+}]=S^{4}_{+}\neq K_{+}$.

Цей приклад показує, що лему 1 та її наслідки не завжди можна
застосовути до гіперкомплексно неопуклих компактів.

Розповсюдимо теорему Клі опуклого аналізу [5] на гіперкомплексний
випадок.

\textbf{Означенняя 4.} \emph{h}\emph{-інтервалом} з центром в точці
\emph{х} радіуса \emph{r} називається перетин відкритої кулі радіуса
\emph{r} з центром в точці \emph{х} з гіперкомплексною прямою, яка
проходить через точку \emph{х}. [7]

\textbf{Означенняя 5.} Точка $x\in E\subset\mathbb{H}^{n}$
називається \emph{h}\emph{-екстремальною точкою множини} \emph{Е},
якщо в \emph{Е} немає жодного \emph{h}-інтервалу, який містить
\emph{х}. [7]

\textbf{Означенняя 6.} \emph{h}\emph{-променем} назвемо замкнену
необмежену ациклічну підмножину гіперкомплексної прямої з
непорожньою межею.

\textbf{Означенняя 7.} \emph{Екстремальним }\emph{h}\emph{-променем}
множини $E\subset\mathbb{H}^{n}$ назвемо \emph{h}-промінь \emph{H},
який належить множині \emph{Е}, якщо множина $E\setminus H$
гіперкомплексно опукла та кожна точка межі променя \emph{H} буде
\emph{h}-екстремальною точкою для множини \emph{Е}. (Це еквівалентне
тому, що жодна точка променя \emph{H} не буде внутрішньою для
довільного \emph{h}-інтервалу, який належить множині \emph{Е} та має
хоча б одну точку за межами \emph{H}.)

Для множини $E\subset\mathbb{H}^{n}$ позначимо: hext~\emph{E}~--
множину її \emph{h}-екстремальних точок, rhext~\emph{E}~-- множину
\emph{h}-екстремальних променів, hconv~\emph{E}~-- \emph{h}-оболонку
\emph{Е}.

\textbf{Лема 2.} Нехай $E\subset\mathbb{H}^{n}$~-- замкнене сильно
гіперкомплексно опукле тіло (тобто $\emph{int}E\neq\emptyset$) з
непорожньою сильно гіперкомплексно опуклою межею $\partial F$, тоді
\emph{Е} має вигляд $E=E_{1}\times\mathbb{H}^{n-1}$, де $E_{1}$~--
ациклічна підмножина прямої $\mathbb{H}$ з непорожньою внутрішністю
відносно цієї прямої.

\emph{Доведення.} Оскільки межа $\partial F$~-- сильно
гіперкомплексно опукла, то для довільної точки $x\in \texttt{int}E$
існує гіперплощина, яка не перетинає $\partial F$. Тому множина
\emph{Е} містить гіперплощину. Отже, за теоремою 3 [3] множину
\emph{Е} можна зобразити у вигляді декартового добутку
$E=E_{1}\times \mathbb{H}^{n-1}$. Множина $E_{1}$ буде ациклічною,
тому що існують перетини \emph{Е} гіперкомплексними прямими, які
гомеоморфні $E_{1}$.

\textbf{Означення 8.} Афінна підмножина \emph{L} називається
\emph{дотичною} до множини \emph{Е}, якщо $L\bigcap
\overline{E}\subset\partial E$, $L\bigcap
\overline{E}\neq\emptyset$.

\textbf{Лема 3.} Якщо $E\subset\mathbb{H}^{n}$~-- сильно
гіперкомплексно опукла замкнена множина та \emph{L}~-- її дотична
гіперкомплексна пряма, то $\texttt{hext}(E\bigcap
L)=(\texttt{hext}E)\bigcap L$.

\emph{Доведення.} Оскільки справедливе включення множин $E\bigcap
L\subset E$, то за означенням \emph{h}-екстремальних точок маємо
$\texttt{hext}(E\bigcap L)\supset(\texttt{hext}E)\bigcap L$. Нехай
точка $x\in\texttt{hext}(E\bigcap L)$. Тоді не може виконуватися
включення $x\in [K]\setminus K$, де $K\subset E$, бо інакше було б
$K\subset E\bigcap L$ (тому що $x\in L$ та \emph{L} є
гіперкомплексною прямою, дотичною до \emph{Е}). А це суперечить
тому, що $x\in\texttt{hext}(E\bigcap L)$. Отже, вірне обернене
включення $\texttt{hext}(E\bigcap L)\subset(\texttt{hext}E)\bigcap
L$ та лема доведена.

\emph{Зауваження 1.} Аналогічно доводиться рівність
$\texttt{rhext}(E\bigcap L)=(\texttt{rhext}E)\bigcap L$ для
\emph{h}-екстремальних променів.

\textbf{Теорема 3.} Кожна замкнена сильно гіперкомплексно опукла
множина $E\subset\mathbb{H}^{n}$, яка не містить гіперкомплекної
прямої, буде \emph{h}-оболонкою своїх \emph{h}-екстремальних точок
та \emph{h}-екстремальних променів $E=\texttt{hconv}(\texttt{hext}
E\bigcup \texttt{rhext} E)$.

\emph{Доведення.} Доведення проведемо індукцією згідно
гіперкомплексної розмірності множини \emph{Е}. При
$\texttt{dim}_{\mathbb{H}}E=0$ та $\texttt{dim}_{\mathbb{H}}E=1$
теорема очевидна. Припустимо, що теорема правильна при всіх
гіперкомплексних розмірностях множини \emph{Е}, які менші \emph{m}
$(1<m\le n)$. Доведемо її при для $\texttt{dim}_{\mathbb{H}}E=m$.

За умовою теореми множина \emph{Е} не містить гіперкомплексної
прямої, тому не може співпадати ні з її афінною оболонкою, ні з
декартовим добутком $E_{1}\times\mathbb{H}^{n-1}$. Тому з леми 2
випливає, що непорожня межа $\partial E$ не буде сильно
гіперкомплексно опуклою множиною.

За означенням сильної гіперкомплексної опуклості перетин множини
\emph{Е} довільною гіперкомплексною прямою буде також сильно
гіперкомплексно опуклим. Нехай \emph{х}~-- довільна точка множини
\emph{Е}. Якщо \emph{х} належить якійсь дотичній прямій \emph{L} до
\emph{Е}, то за припущенням індукції маємо включення $x\in
\texttt{hconv}((\texttt{hext}E\bigcap L)\bigcup
\texttt{rhext}(E\bigcap L))$. Якщо ж існують точки множини \emph{Е},
через які не проходить жодна дотична до \emph{Е} гіперкомплексна
пряма, тоді знайдеться точка $x\in \texttt{int}E$.

У цьому випадку через точку \emph{х} проведемо гіперкомплексну пряму
\emph{l}. Перетин $l\bigcap E$ є сильно гіперкомплексно опуклою
множиною та не співпадає з \emph{l}. Тому $x\notin
[\partial(l\bigcap E)]$. Нехай тепер \emph{y}~-- довільна точка межі
перетину $\partial(l\bigcap E)$. Згідно сильної гіперкомплексної
опуклості, через точку \emph{y} можна провести дотичну до множини
\emph{Е} пряму \emph{Т}. За припущенням індукції отримаємо включення
$y\in \texttt{hconv}((\texttt{hext}E\bigcap T)\bigcup
\texttt{rhext}(E\bigcap T))$. Зауважимо, що це виконується для
кожної точки $y\in\partial(l\bigcap E)$. Тоді, враховуючи лему 3 та
зауваження 2, отримаємо $x\in \texttt{hconv}(\texttt{hext}E\bigcup
\texttt{rhext}E)$. В силу довільності вибору точки \emph{х}
$E\subset \texttt{hconv}(\texttt{hext}E\bigcup \texttt{rhext}E)$.
Обернене включення тривіальне. Теорема доведена.

Клас сильно гіперкомплексно опуклих множин є незамкненим відносно
перетинів. Тому для нього не виконується основна аксіома
опуклості~-- перетин будь-якої кількості опуклих множин повинен бути
опуклим. Означимо клас множин, який включає в себе сильно
гіперкомплексно опуклі множини та є замкненим відносно перетинів.

\textbf{Означення 9.} Гіперкомплексно опуклу множину
$E\subset\mathbb{H}^{n}$ назвемо $\mathbb{H}$-\emph{квазіопуклою
множиною}, якщо її перетин довільною гіперкомплексною прямою
$\gamma$ не містить тривимірного коциклу, тобто $H^{3}(\gamma\bigcap
E)=0$.

Очевидно, що клас $\mathbb{H}$-квазіопуклих множин включає в себе
сильно гіперкомплексно опуклі області та компакти.

Покажемо замкненість класу $\mathbb{H}$-квазіопуклих множин в тому
сенсі, що перетин довільної сім'ї компактних
$\mathbb{H}$-квазіопуклих множин буде $\mathbb{H}$-квазіопуклою
множиною.

\textbf{Теорема 4.} Перетин довільної сім'ї
$\mathbb{H}$-квазіопуклих компактів буде $\mathbb{H}$-квазіопуклим
компактом.

\emph{Доведення.} Доведення достатньо провести для двох компактів.
Нехай $K_{1},K_{2}$~-- два довільні $\mathbb{H}$-квазіопуклі
компакти, $\gamma$~-- довільна гіперкомплексна пряма, яка перетинає
множину $K_{1}\bigcap K_{2}$. Використаємо точну когомологічну
послідовність Майєра-В'єторіса
$$H^{3}(\gamma\bigcap K_{1})\bigoplus H^{3}(\gamma\bigcap
K_{2})\rightarrow H^{3}(\gamma\bigcap K_{1}\bigcap K_{2})\rightarrow
H^{4}(\gamma\bigcap (K_{1}\bigcup K_{2})).$$

Оскільки компакти $K_{1}$ та $K_{2}$ $\mathbb{H}$-квазіопуклі, то
$H^{3}(\gamma\bigcap K_{1})=0$ та $H^{3}(\gamma\bigcap K_{2})=0$.
Тому $H^{3}(\gamma\bigcap K_{1})\bigoplus H^{3}(\gamma\bigcap
K_{2})=0$.

З іншого боку, компактний перетин $\gamma\bigcap (K_{1}\bigcup
K_{2})=(\gamma\bigcap K_{1})\bigcup(\gamma\bigcap K_{2})$ не може
містити всю гіперкомплексну пряму $\gamma$, яка є чотиривимірним
дійсним многовидом, тому $H^{4}(\gamma\bigcap (K_{1}\bigcup
K_{2}))=0$.

З точності когомологічної послідовності випливає, що
$H^{3}(\gamma\bigcap K_{1}\bigcap K_{2})=0$. Це еквівалентне тому,
що перетин множини $K_{1}\bigcap K_{2}$ довільною гіперкомплексною
прямою не містить тривимірного коциклу. Звідси випливає
$\mathbb{H}$-квазіопуклість компакта $K_{1}\bigcap K_{2}$. Теорема
доведена.

\textbf{Лема 4.} Якщо всі співмножники $E_{j}, j=1,...,n,$
ациклічні, то довільний перетин множини $E=E_{1}\times ...\times
E_{n}\subset\mathbb{H}^{n}$ не містить тривимірного коциклу.

\emph{Доведення.} Доведення проведемо для випадку, коли
$E=E_{1}\times E_{2}\subset\mathbb{H}^{2}$ та гіперкомплексна пряма
$\gamma$ проходить через початок координат. Цього завжди можна
досягти заміною координат. Рівняння прямої має вигляд
$$\gamma=\{x|x_{1}=k_{1}x_{2},x_{3}=...=x_{n}=0,k=-\frac{r_{1}}{r_{2}},r_{i}=|x_{i}|,i=1,2\}.$$
Оскільки в прямій $\gamma$ координати $x_{1},x_{2}$ пов'язані між
собою співвідношеннями $x_{1}=kx_{2}$, а множини $E_{2}$ та
$kE_{2},k\neq0$ гомеоморфні між собою, то перетин $\gamma\bigcap E$
співпадає з $E_{1}\bigcap kE_{2}$.

З точності когомологічної послідовності Майєра-В'єторіса
$$H^{3}(E_{1})\bigoplus H^{3}(kE_{2})\rightarrow H^{3}(E_{1}\bigcap kE_{2})\rightarrow
H^{4}(E_{1}\bigcup kE_{2}),$$ оскільки $H^{3}(E_{1})\bigoplus
H^{3}(kE_{2})=0$ та $H^{4}(E_{1}\bigcup kE_{2})=0$, отримуємо, що
$H^{3}(E_{1}\bigcap kE_{2})=0$. Лема доведена.

З леми 4 фактично випливає, що декартовий добуток
$\mathbb{H}$-квазіопуклих компактів буде $\mathbb{H}$-квазіопуклою
множиною.

\textbf{Означення 10.} \emph{Лінійним поліедром} називається множина
виду $E=\{x|f_{j}(x)\in E_{j},j\in J=\{1,2,...,N\}\}$, де
$E_{j}\subset\mathbb{H}^{1}$, $f_{j}(x)=\sum_{k=1}^{n}a_{jk}x_{k}$,
причому довільні дві функції $f_{k}(x),f_{j}(x),k\neq j,$ є лінійно
незалежними, а кожна з функцій $f_{j}$ відображає \emph{Е} в
підмножину гіперкомплексної прямої $E_{j}$. [3]

Оскільки довільний компактний лінійний поліедр можна зобразити у
вигляді перетину не більше, ніж \emph{n}, (\emph{n}~-- скінченне
число і позначає кількість граней поліедра) декартових добутків
$\Gamma_{i}\times K_{i}$, $i=\overline{1,n}$, множини $\Gamma_{i}$
де $\Gamma_{i}$~-- грань поліедра, на кулю $K_{i}$ досить великого
радіуса в $(n-1)$-вимірній гіперкомплексній гіперплощині, яка
ортогональна до грані $\Gamma_{i}$, то з теореми 4 отримуємо
наступне твердження.

\textbf{Теорема 5.} Компактний лінійний поліедр, всі грані якого не
містять трьохвимірних циклів, є $\mathbb{H}$-квазіопуклою множиною.

\textbf{Наслідок 5.} Перетин сильно гіперкомплексно опуклих
компактів буде $\mathbb{H}$-квазіопуклою множиною.

\emph{Зауваження 2.} У роботі [8] доведено, що гіперкомплексно
опукла область з гладкою межею буде $\mathbb{H}$-квазіопуклою
множиною.

\textbf{Теорема 6.} Кожен ациклічний в розмірності три
гіперкомплексно опуклий компакт \emph{Е} буде
$\mathbb{H}$-квазіопуклим.

\emph{Доведення.} Розглянемо перетин множини \emph{Е}
гіперкомплексною прямою $\gamma$. Нехай перетин $\gamma\bigcap E$
містить негомологічний нулю тривимірний цикл \emph{Z}. Візьмемо
точку $x_{0}\in\gamma\setminus\gamma\bigcap E$, зачеплену з циклом
\emph{Z}. Оскільки множина \emph{Е} гіперкомплексно опукла, то існує
гіперкомплексна гіперплощина \emph{L}, яка проходить через точку
$x_{0}$ і не перетинає \emph{Е}. Тоді $(4n-4)$-вимірний цикл
$K=L\bigcup\infty$ зачеплений з трьохвимірним циклом \emph{Z}.
Оскільки множина \emph{Е} однозв'язна, то цикл \emph{Z} гомологічний
нулю в \emph{Е}. Тоді існує ланцюг \emph{С} в \emph{Е}, який
обмежується циклом \emph{Z}. За теоремою про зачеплення $C\bigcap
K\neq\emptyset$, що суперечить тому, що множина $K\bigcap E$
порожня. Теорема доведена.


\vspace{10mm}

\noindent {\bf М.~В.~Стефанчук }  \vspace{3mm} \\
\noindent {\bf Экстремальные элементы в гиперкомплексном пространстве }\\

\noindent В работе исследуются экстремальные элементы и
\emph{h}-оболочка множеств в \emph{n}-мерном гиперкомплексном
пространстве $\mathbb{H^{\emph{n}}}$. Вводится класс
$\mathbb{H}$-квазивыпуклых множеств, которые включают в себя сильно
гиперкомплексно выпуклые множества и являются замкнутыми
относительно пересечений.

\textbf{Ключевые слова:} гиперкомплексно выпуклое множество, сильно
гиперкомплексно выпуклое множество, \emph{h}-оболочка множества,
\emph{h}-экстремальная точка, \emph{h}-экстремальный луч,
$\mathbb{H}$-квазивыпуклое множество.

\vspace{10mm}

\noindent {\bf M.~V.~Stefanchuk }  \vspace{3mm} \\
\noindent {\bf Extremal elements in hypercomplex
space}\\

\noindent Extremal elements and a \emph{h}-hull of sets in the \emph{n}-dimensional hypercomplex
space $\mathbb{H^{\emph{n}}}$ are investigated. Introduced  a class of  $\mathbb{H}$-quasiconvex
sets including strongly hypercomplex convex sets and being closed with respect to intersections.


\vspace{2mm}

\noindent
Інститут математики НАН України, Київ\\
Адреса: 01601, Київ, вул. Терещенківська, 3.\\
Автор:  Стефанчук Марія Володимирівна\\
e-mail: stefanmv43@gmail.com

\end{document}